\documentclass[12pt]{amsart}
\usepackage{amssymb}
\makeatletter
\renewcommand{\section}{\@startsection{section}{1}%
  \z@{.7\linespacing\@plus\linespacing}%
  {-.5em}%                     % si < 0 pas de passage a la ligne
  {\bfseries\large}}           % Le style (taille, bold ...)

\renewcommand{\subsection}{\@startsection{subsection}{2}%
  \z@{.5\linespacing\@plus.7\linespacing}{-.5em}%
  {\normalfont\bfseries}}

\renewcommand{\subsubsection}{\@startsection{subsubsection}{3}%
  \z@{.5\linespacing\@plus.7\linespacing}{-.5em}%
  {\normalfont\itshape}}
\makeatother

\newenvironment{Prf}{\vskip 1em{\it Proof} :}%
{\unskip\hfill\null\nobreak\hfill\carre\vskip1em\par}
\newcommand{\carre}{\rule{1ex}{1ex}}

\def \rb {\Bbb{R}}
\def \nb {\Bbb{N}}
\def \zb {\Bbb{Z}}

\begin{document}

{\bf Subharmonic solutions for nonautonomous
 sublinear first order\\ Hamiltonian systems}\\
 
 \bigskip

{\bf A. Raouf Chouikha }\footnote { Universite Paris 13 LAGA, Villetaneuse 93430, chouikha@math.univ-paris13.fr} and  {\bf Mohsen Timoumi}\footnote { Faculte des Sciences de Monastir, Tunisie, m\_timoumi@yahoo.com.}\\
\bigskip

{\bf Abstract}\quad  
In this paper, the existence of subharmonic solutions for a class of non-autonomous first-order Hamiltonian systems is investigated. We also study the minimality of periods for such solutions. Our results which extend and improve many previous results will be illustrated by specific examples. Our main tools are the minimax methods in critical point theory and the least action principle.\\
{\bf Key words.} Hamiltonian systems. Critical point theory. Least action principle. Subharmonic solutions.\\
\bigskip

\section{Introduction} Consider the nonautonomous first-order Hamiltonian system
$$\dot x(t)= JH'(t,x(t)) \leqno (\mathcal H)$$
where $H: \rb\times\rb^{2N}\longrightarrow\rb,\ (t,x)\longmapsto H(t,x)$ is a continuous function, $T-$ periodic $(T > 0)$ in the first variable and differentiable with respect to the second variable with continuous derivative $H'(t,x)={{\partial H}\over{\partial x}}(t,x)$ and $J$ is the standard symplectic matrix:
$$J=\left(\begin{array}{ll}
0 & -I_{N}\\
I_{N} & 0
\end{array}\right)$$
$I_N$ being the identity matrix of order $N$.\\
In this work, we are focused in the existence of subharmonic solutions of $(\mathcal H)$. Assuming that $T>0$ is the minimal period of the time dependence of $H(t,x)$, by subharmonic solution of $(\mathcal H)$ we mean a $kT-$periodic solution, where $k$ is any integer; when moreover the periodic solution is not $T-$periodic we call it a true subharmonic.\\
Considerable attention has been paid in the last years to the subharmonic solutions of Hamiltonian systems. Most research on subharmonics concern second order systems. Indeed, several papers have been published in this direction, we refer the reader to [3,4,8,14,18] and references therein. Concerning the first order, few researchers are interested because the problem difficult at first. Note, however, the following works [1,2,6,7,9]. Using variational methods, many papers devoted to the existence of subharmonics for $(\mathcal H)$ with various assumptions on the growth of the Hamiltonian. In particular, under the assumptions that there exists a constant $M > 0$ such that
$$\left|H'(t,x)\right|\leq M,\ \forall x\in\rb^{2N},\ \forall t\in[0,T],\leqno(1.1)$$
and
$$\lim_{\left|x\right|\longrightarrow\infty}H(t,x)=\pm\infty,\ uniformly\ in\ t\in [0,T],\leqno(1.2)$$
[16] has shown that the system $(\mathcal H)$ admitted a sequence of subharmonic solutions. After that, [1] generalized this result to the sublinear case. Precisely, it was assumed that the nonlinearity satisfied the following restrictions:\\
$$\left|H'(t,x)\right|\leq f(t)\left|x\right|^{\alpha}+g(t),\ \forall x\in\rb^{2N},\ a.e.\ t\in[0,T], \leqno(1.3)$$
$$\frac{1}{\left|x\right|^{2\alpha}}\int^{T}_{0}H(t,x)dt\longrightarrow+\infty\ as\ \left|x\right|\longrightarrow +\infty, \leqno(1.4)$$
where $f\in L^{\frac{2}{1-\alpha}}(0,T;\rb^{+})$ and $g\in L^{2}(0,T;\rb^{+})$ are $T-$ periodic and $\alpha\in[0,1[$. \\
(1.5) There exists a subset $C$ of $[0,T]$ with $meas(C)>0$ and a $T-$ periodic function $f\in L^{1}(0,T;\rb)$ such that
$$H(t,x)\longrightarrow+\infty\ as\ \left|x\right|\longrightarrow\infty,\ for\ a.e.\ t\in C,$$
and
$$H(t,x)\geq f(t),\ \forall x\in\rb^{2N},\ a.e.\ t\in [0,T].$$
Under these conditions, subharmonic solutions of the system $(\mathcal H)$ have been obtained. More precisely, it was proved that for all integer $k\geq 1$, the system $(\mathcal H)$ possesses a $kT-$ periodic solution $x_{k}$ such that $\left\|x_{k}\right\|_{\infty}\longrightarrow\infty$ as $k\longrightarrow\infty$.\\

Our paper is organized as follows. In section 3, we will be interested in the existence of subharmonics of $(\mathcal H)$ under some more general conditions than (1.3), (1.4). In section 4, we will study the minimality of periods of the subharmonic solutions. We will give examples in order to show the originality of our results which improve many previous results among them [1,14,16].\\
For the proofs, we will apply a Generalized Saddle Point Theorem to the Least Action Integral and use a Generalized Egoroff's Lemma.\\
\section{Preliminaries} Firstly, let us recall a critical point theorem due to [5] which will be useful in the sequel.\\
Let $E=W\oplus Z$ be a Banach space and $(E_{n}=W_n\oplus Z_n)$ be a sequence of closed subspaces with 
$Z_{0}\subset Z_{1}\subset...\subset Z$, $W_{0}\subset W_{1}\subset...\subset W$, $1 \leq dim\ W_{n} < +\infty$. For $f\in C^{1}(E,\rb)$, we denote by $f_{n} = f_{|E_{n}}$ the restriction of $f$ into $E_{n}$. Then we have $f_{n} \in C^1 (E_{n},\rb)$, for all $n\geq 1$.\\
{\bf Definition 2.1.} Let $f \in C^{1} (E,\rb)$ and $c\in\rb$. The function $f$ satisfies the Palais-Smale condition with respect to $(X_n)$ at a level $c \in \rb$ if every sequence $(x_n)$ satisfying
$$n_{j}\longrightarrow\infty,\ x_{n_{j}}\in E_{n_{j}},\ f(x_{n_{j}}) \longrightarrow c,\ f'_{n_{j}} (x_{n_{j}})\longrightarrow 0$$
possesses a subsequence which converges in $E$ to a critical point of $f$. The above property will be referred as the $(PS)^{*}_{c}$ condition with respect to $(E_{n})$.\\
{\bf Theorem 2.1} (Generalized Saddle Point Theorem). Let $f\in C^{1}(E,\rb)$. Assume that there exists a constant $r > 0$ such that with $Y=\left\{w\in W:\left\|w\right\|=r\right\}$:\\
$$\sup_{Y}f\leq \inf_{Z}f, \leqno a)$$
b) $f$ is bounded above on $A=\left\{w\in W:\left\|w\right\|\leq r\right\},$\\
c) $f$ satisfies the $(PS)_c^*$, with
$$c=\inf_{A\in\mathcal A}\sup_{x\in A}f(x),$$
with
$$\mathcal{A}=\left\{A\subset E:A\ is\ closed,\ Y\subset A,\ cat_{E,Y}(A)=1\right\}.$$
Then $c$ is a critical value of $f$ and $c\geq \inf_{Z}f$.\\
{\bf Remark 2.1.} In a) we may replace $Z$ by $q+Z$, $q\in W$.\\
Consider the Hilbert space $E = H^{1\over 2}(S^1,\rb^{2N})$ where $S^1 = {\rb/{(T\zb)}}$ and the continuous quadratic form $Q$ defined in $E$ by
$$Q(u)=\frac{1}{2}\int^{T}_{0}J\dot{u}\cdot u dt$$
where $x\cdot y$ inside the sign integral is the inner product of $x,y\in\rb^{2N}$. Let us denote by $E^{0}$, $E^{-}$, $E^{+}$ respectively the subspaces of $E$ on which $Q$ is null, negative definite and positive definite. It is well known that these subspaces are mutually orthogonal in $L^2(S^1,\rb^{2N})$ and in $E$ with respect to the bilinear form:
$$B(u,v)=\frac{1}{2}\int^{T}_{0}J\dot{u}\cdot v dt,\ u,v\in E$$
associated to $Q$. If $u\in E^+$ and $v\in E^-$, then $B(u,v)=0$ and $Q(u+v)=Q(u)+Q(v)$. For $u=u^{+} +u^{-} +u^{0}\in E$, the expression
$$||u||=\big[Q(u^{+})-Q(u^{-})+|u^{0}|^{2}\big]^{\frac{1}{2}}$$
is an equivalent norm in $E$. Moreover, the space $E$ is compactly embedded in $L^{s}(S^1,\rb^{2N})$ for all $s\in [1,\infty[$ (see [11]). In particular for all $s\in [1,\infty[$, there exists a constant $\lambda_{s}>0$ such that for all $u\in E$,
$$||u||_{L^{s}} \leq \lambda_{s}||u||.\leqno(2.1)$$
\section{Existence of subharmonics} Let $\gamma:\rb^{+}\longrightarrow\rb^{+}$ be a nondecreasing continuous function satisfying the properties:
$$\gamma(s+t)\leq c(\gamma(s)+\gamma(t)),\ \forall s,t\in\rb^{+},\leqno(i)$$
$$0\leq\gamma(t)\leq at^{\alpha}+b,\ \forall t\in\rb^{+},\leqno(ii)$$
$$\gamma(t)\longrightarrow +\infty\ as\ t\longrightarrow+\infty,\leqno(iii)$$
where $a,b,c$ are positive constants and $\alpha\in [0,1[$. Consider the following assumptions:\\
$(H_{1})$ There exist two $T-$periodic functions $p\in L^{2\over{1-\alpha}}(0,T;\rb^{+})$ and $q\in L^{2}(0,T;\rb^{+})$ such that \\
$$\left|H'(t,x)\right|\leq p(t)\gamma(\left|x\right|)+q(t),\ \forall x\in\rb^{2N},\ a.e.\ t\in [0,T];$$
$(H_{2})$ Either
$$\frac{1}{\gamma^{2}(\left|x\right|)}\int^{T}_{0}H(t,x)dt\longrightarrow+\infty\ as\ \left|x\right|\longrightarrow +\infty, \leqno(i)$$
or
$$\frac{1}{\gamma^{2}(\left|x\right|)}\int^{T}_{0}H(t,x)dt\longrightarrow-\infty\ as\ \left|x\right|\longrightarrow +\infty; \leqno(ii)$$
$(H_{3})$ There exist a subset $C$ of $[0,T]$ with $meas(C)>0$ and a $T-$ periodic function $f\in L^{1}(0,T;\rb)$ such that either
$$H(t,x)\longrightarrow+\infty\ as\ \left|x\right|\longrightarrow\infty,\ a.e.\ t\in C, \leqno(i)$$
and
$$H(t,x)\geq f(t),\ \forall x\in\rb^{2N},\ a.e.\ t\in[0,T];$$
or
$$H(t,x)\longrightarrow-\infty\ as\ \left|x\right|\longrightarrow\infty,\ a.e.\ t\in C, \leqno(ii)$$
and
$$H(t,x)\leq f(t),\ \forall x\in\rb^{2N},\ a.e.\ t\in[0,T].$$
Our main result in this section reads as follows.\\
{\bf Theorem 3.1.} Suppose $(H_{1})-(H_{3})$ hold. Then, for all positive integer $k$, the Hamiltonian system $(\mathcal H)$ possesses at least one $kT-$ periodic solution $x_{k}$ satisfying
$$\lim_{k\longrightarrow\infty}\left\|x_{k}\right\|_{\infty}=+\infty$$
where $\left\|x\right\|_{\infty}=\sup_{t\in\rb}\left|x(t)\right|$.\\
{\bf Corollary 3.1.} Assume $H$ satisfies assumption $(H_{1})$ and\\
$(H_{4})$ There exist a subset $C$ of $[0,T]$ with $meas(C)>0$ and a $T-$ periodic function $f\in L^{1}(0,T;\rb)$ such that either
$$\frac{H(t,x)}{\gamma^{2}(\left|x\right|)}\longrightarrow+\infty\ as\ \left|x\right|\longrightarrow\infty,\ a.e.\ t\in C,\leqno(i)$$
and
$$H(t,x)\geq f(t),\ \forall x\in\rb^{2N},\ a.e.\ t\in[0,T],$$
or
$$\frac{H(t,x)}{\gamma^{2}(\left|x\right|)}\longrightarrow-\infty\ as\ \left|x\right|\longrightarrow\infty,\ a.e.\ t\in C,\leqno(ii)$$
and
$$H(t,x)\leq f(t),\ \forall x\in\rb^{2N},\ a.e.\ t\in[0,T].$$
Then the conclusion of Theorem 3.1 holds.\\
{\bf Example 3.1.} Theorem 3.1 in [16] and Theorem 1.1 in [1] are special cases of Theorem 3.1 with control function $\gamma(t)=t^{\alpha}$, $0\leq \alpha<1$, $t\in\rb^{+}$. What's more, there are functions $H(t,x)$ satisfying our theorem and do not satisfy the results in [1,16]. For example, we consider the Hamiltonian
$$H(t,x)=\theta(t)ln^{\frac{3}{2}}(1+\left|x\right|^{2}),$$
where $\theta$ is the $T-$ periodic function such that its restriction to $[0,T]$ is given by.
$$\theta(t)=\left\{
\begin{array}{l}
sin(\frac{2\pi}{T}t),\ t\in [0,\frac{T}{2}]\\
0,\ \ \ \ \ \ \ \ \ \ \ t\in[\frac{T}{2},T],
\end{array}\right.$$
It is clear that $H(t,x)$ does not satisfy (1.1), (1.2) nor (1.3), (1.4), (1.5). Take $\gamma(t)=ln^{\frac{1}{2}}(1+t^{2})$. It is not difficult to see that $\gamma$ is nondecreasing and satisfies $(i)$ and $(iii)$. For $(ii)$, we have\\
$$1+(s+t)^{2}\leq (1+s^{2})^{2}(1+t^{2})^{2},$$
and since $ln$ is increasing, we get
$$ln(1+(s+t)^{2})\leq 2[ln(1+s^{2})+ln(1+t^{2})]$$
which with the property $\sqrt{a+b}\leq \sqrt{a}+\sqrt{b}, \forall a,b\in\rb^{+}$ yield
$$ln^{\frac{1}{2}}(1+(s+t)^{2})\leq \sqrt{2}[ln^{\frac{1}{2}}(1+s^{2})+ln^{\frac{1}{2}}(1+t^{2})].$$
It is easy to verify that $H$ satisfies $(H_{1})$, $(H_{2})$ and $(H_{3})$ with $C=]0,\frac{T}{2}[$.\\
{\bf Proof of Theorem 3.1.} Firstly, let us remark the following:\\
{\bf Remark 3.1.} Let $x(t)$ be a periodic solution of $(\mathcal H)$, then by replacing $t$ by $-t$ in $(\mathcal H)$, we obtain
$$\dot{x}(-t) = JH'(-t,x(-t)).$$
So it is clear that the function $y(t)=x(-t)$ is a periodic solution of the system
$$\dot{y}(t) = -JH'(-t,y(t)).$$
Moreover, $-H(-t,x)$ satisfies $(H_{2})(i)-(H_{4})(i)$ whenever $H(t,x)$ satisfies respectively $(H_{2})(ii)-(H_{4})(ii)$. Hence, in the following, we will assume that $H$ satisfies $(H_{1})$, $(H_{2})(i)$ and $(H_{3})(i)$.\\
By making the change of variables $t\longrightarrow\frac{t}{k}$, the system $(\mathcal H)$ transforms to
$$\dot{u}(t)=kJH'(kt,u(t)). \leqno({\mathcal H}_{k})$$
Hence, to find $kT-$ periodic solutions of $(\mathcal H)$, it suffices to find $T-$ periodic solutions of $({\mathcal H}_{k})$.\\
Consider the family of functionals $(\Phi_{k})_{k\in\nb}$ defined on the space $E$ introduced above by
$$\Phi_{k}(u)=\int^{T}_{0}[\frac{1}{2}J\dot{u}(t)\cdot u(t)+kH(kt,u(t))]dt.$$
By assumption $(H_{1})$ and the property $(ii)$ of $\gamma$, we have
$$\left|H'(t,x)\right|\leq p(t)[a\left|x\right|^{\alpha}+b]+q(t),\ \forall x\in\rb^{2N},\ a.e.\ t\in[0,T]. \leqno(3.1)$$
So, by Proposition B.37 [11], $\Phi_{k}\in C^{1}(E,\rb)$ and critical points of $\Phi_{k}$ on $E$ correspond to the $T-$ periodic solutions of $({\mathcal H}_{k})$, moreover one has
$$\Phi'_{k}(u)v=\int^{T}_{0}[J\dot{u}(t)+kH'(kt,u(t))]\cdot v(t)dt,\ \forall u,v\in E. \leqno(3.2)$$
Let us fix a positive integer $k$, we will study the existence of critical points of the functional $\Phi_{k}$. To this aim we will apply the Generalized Saddle Point Theorem to the functional $\Phi_{k}$ with the decomposition $W=E^{-}$, $Z=E^{0}\oplus E^{+}$ of $E$ and with respect to the sequence of subspaces
$$E_{n}=\{u\in E/\ u(t)=\sum_{\left|m\right|\leq n}exp({{2\pi}\over T}mtJ)\hat{u}_{m}\ a.e.\},\ n\geq 0.$$
Firstly, let us check the Palais-Smale condition.\\
{\bf Lemma 3.1.} For all level $c\in\rb$, the functional $\Phi_{k}$ satisfies the $(PS)^{*}_{c}$ condition with respect to the sequence $(E_{n})_{n\in\nb}$.\\
\begin{Prf}
Let $c\in\rb$ and let $(u_{n})_{n\in\nb}$ be a sequence such that for a subsequence $(n_{j})$ of $\nb$
$$n_{j}\longrightarrow\infty,\ u_{n_{j}}\in E_{n_{j}},\ \Phi_{k}(u_{n_{j}})\longrightarrow c\ and\ \Phi'_{k,n_{j}}(u_{n_{j}}) \longrightarrow 0\ as\ j\longrightarrow\infty,$$
where $\Phi_{k,n_{j}}$ is the functional $\Phi_{k}$ restricted to $E_{n_{j}}$. Set $u_{n_{j}}=\bar{u}_{n_{j}}+\tilde{u}_{n_{j}}$, with $\bar{u}_{n_{j}}=\frac{1}{T}\int^{T}_{0}u_{n_{j}}(t)dt$ and $\tilde{u}_{n_{j}}=u_{n_{j}}-\bar{u}_{n_{j}}$, we have the relation
$$\Phi'_{k,n_{j}}(u_{n_{j}})(u^{+}_{n_{j}}-u^{-}_{n_{j}})=2\left\|\tilde{u}_{n_{j}}\right\|^{2}+k\int^{T}_{0}H'(kt,u_{n_{j}}).(u^{+}_{n_{j}}-u^{-}_{n_{j}})dt.$$
Since $\Phi'_{k,n_{j}}(u_{n_{j}})\longrightarrow 0$ as $j\longrightarrow\infty$, there exists a constant $c_{1}>0$ such that
$$\left|\Phi'_{k,n_{j}}(u_{n_{j}})(u^{+}_{n_{j}}-u^{-}_{n_{j}})\right|\leq c_{1}\left\|\tilde{u}_{n_{j}}\right\|,\ \forall j\in\nb. \leqno(3.3)$$
By H$\ddot{o}$lder's inequality and $(H_{1})$
$$\left|\int^{T}_{0}H'(kt,u_{n_{j}}).-(u^{+}_{n_{j}}-u^{-}_{n_{j}})dt\right|\leq                               \left\|\tilde{u}_{n_{j}}\right\|_{L^{2}}(\int^{T}_{0}\left|H'(kt,u_{n_{j}})\right|^{2}dt)^{\frac{1}{2}}$$
$$\leq\left\|\tilde{u}_{n_{j}}\right\|_{L^{2}}(\int^{T}_{0}[p(kt)\gamma(\left|u_{n_{j}}\right|)+q(kt)]dt)^{\frac{1}{2}}$$ 
$$\leq\left\|\tilde{u}_{n_{j}}\right\|_{L^{2}}\big[(\int^{T}_{0}p^{2}(kt)\gamma^{2}(\left|u_{n_{j}}\right|)dt)^{\frac{1}{2}}+\left\|q\right\|_{L^{2}}\big]. \leqno(3.4)$$ 
Now, by the nondecreasing and the properties (i) and (ii) of $\gamma$, we have
$$(\int^{T}_{0}p^{2}(kt)\gamma^{2}(\left|u_{n_{j}}\right|)dt)^{\frac{1}{2}}\leq(\int^{T}_{0}p^{2}(kt)\gamma^{2}(\left|\tilde{u}_{n_{j}}\right|+\left|\bar{u}_{n_{j}}\right|)dt)^{\frac{1}{2}}$$
$$\leq c(\int^{T}_{0} [p^{2}(kt)[\gamma(\left|\tilde{u}_{n_{j}}\right|) +\gamma(\left|\bar{u}_{n_{j}}\right|)]^{2}dt)^{\frac{1}{2}}$$
$$\leq c\big[(\int^{T}_{0}p^{2}(kt)\gamma^{2}(\left|\tilde{u}_{n_{j}}\right|)dt)^{\frac{1}{2}} +\left\|p\right\|_{L^{2}}\gamma(\left|\bar{u}_{n_{j}}\right|)\big]$$
$$\leq c\big[\big(\int^{T}_{0}p^{2}(kt)(a\left|\tilde{u}_{n_{j}}\right|^{\alpha}+b)^{2}dt\big)^{\frac{1}{2}}+ \left\|p\right\|_{L^{2}} \gamma(\left|\bar{u}_{n_{j}}\right|)\big]$$
$$\leq c\big[a\big(\int^{T}_{0}p^{2}(kt)\left|\tilde{u}_{n_{j}}\right|^{2\alpha}dt\big)^{\frac{1}{2}}+ b\left\|p\right\|_{L^{2}}+ \left\|p\right\|_{L^{2}} \gamma(\left|\bar{u}_{n_{j}}\right|)\big]$$
$$\leq c\big[a\left\|p\right\|_{L^{\frac{2}{1-\alpha}}}\left\|\tilde{u}_{n_{j}}\right\|^{\alpha}_{L^{2}}+ b\left\|p\right\|_{L^{2}}+\left\|p\right\|_{L^{2}}\gamma(\left|\bar{u}_{n_{j}}\right|)\big]. \leqno(3.5)$$
Therefore by (2.1), (3.4) and (3.5), there exists a positive constant $c_{2}$ such that
$$k\left|\int^{T}_{0}H'(kt,u_{n_{j}}).(u^{+}_{n_{j}}-u^{-}_{n_{j}})dt\right|\leq c_{2}\left\|\tilde{u}_{n_{j}}\right\|\big[ \left\|\tilde{u}_{n_{j}}\right\|^{\alpha}+ \gamma(\left|\bar{u}_{n_{j}}\right|)+1\big]$$
which with (3.3) yield
$$c_{1}\geq 2\left\|\tilde{u}_{n_{j}}\right\|-c_{2}\big[\left\|\tilde{u}_{n_{j}}\right\|^{\alpha}+ \gamma(\left|\bar{u}_{n_{j}}\right|)+1\big]$$
and
$$c_{2}\gamma(\left|\bar{u}_{n_{j}}\right|)\geq \left\|\tilde{u}_{n_{j}}\right\|\big[2-c_{2}\left\|\tilde{u}_{n_{j}}\right\|^{\alpha-1}\big]-c_{1}-c_{2}. \leqno(3.6)$$
Assume that $(\tilde{u}_{n_{j}})$ is unbounded, then by going to a subsequence, if necessary, we can assume that $\left\|\tilde{u}_{n_{j}}\right\|\longrightarrow\infty$ as $j\longrightarrow\infty$. Since $0\leq\alpha<1$, we deduce from (3.6) that there exists constant $c_{3}>0$ such that
$$\left\|\tilde{u}_{n_{j}}\right\|\leq c_{3}(\gamma(\left|\bar{u}_{n_{j}}\right|)+1) \leqno(3.7)$$
for $j$ large enough. By the continuity of $\gamma$ and $(3.7)$, we have $\left|\bar{u}_{n_{j}}\right|\longrightarrow\infty$ as $j\longrightarrow\infty$.\\
Now, by the Mean Value Theorem, H$\ddot{o}$lder's inequality, properties (2.1), (3.7), property $(ii)$ of $\gamma$ and since $\alpha<1$, we obtain as above
$$k\left|\int^{T}_{0}(H(kt,u_{n_{j}})-H(kt,\bar{u}_{n_{j}}))dt\right|=k\left|\int^{T}_{0}\int^{1}_{0}H'(kt,\bar{u}_{n_{j}}+s\tilde{u}_{n_{j}})\cdot.\tilde{u}_{n_{j}}dsdt\right|$$
$$k\leq\left\|\tilde{u}_{n_{j}}\right\|_{L^{2}}\int^{1}_{0}\big(\int^{T}_{0}\left| H'(kt,\bar{u}_{n_{j}}+s\tilde{u}_{n_{j}})\right|^{2}dt\big)^{\frac{1}{2}}ds$$
$$\leq c_{4}\left\|\tilde{u}_{n_{j}}\right\|\big[ \left\|\tilde{u}_{n_{j}}\right\|^{\alpha}+ \gamma(\left|\bar{u}_{n_{j}}\right|)+1\big]$$
$$\leq c_{4}(\gamma(\left|\bar{u}_{n_{j}}\right|)+1)\big[c^{\alpha}_{3}(\gamma(\left|\bar{u}_{n_{j}}\right|)+1)^{\alpha}+\gamma(\left|\bar{u}_{n_{j}}\right|)+1\big]$$
$$\leq c_{5}\gamma^{2}((\left|\bar{u}_{n_{j}}\right|)\leqno(3.8)$$
for $j$ large enough, where $c_{4}$, $c_{5}$ are two positive constants. Therefore by (3.7), (3.8) there exists a positive constant $c_{6}$ such that for $j$ large enough
$$\Phi_{k}(u_{n_{j}})=\left\|u^{+}_{n_{j}}\right\|^{2}-\left\|u^{-}_{n_{j}}\right\|^{2}+k\int^{T}_{0}[H(kt,u_{n_{j}})-H(kt,\bar{u}_{n_{j}})]dt+k\int^{T}_{0}H(kt,\bar{u}_{n_{j}})dt$$
$$\geq -c^{2}_{3}(\gamma(\left|\bar{u}_{n_{j}}\right|)+1)^{1}-c_{5}\gamma^{2}(\left|\bar{u}_{n_{j}}\right|)+k\int^{T}_{0}H(kt,\bar{u}_{n_{j}})dt$$
$$=\gamma^{2}(\left|\bar{u}_{n_{j}}\right|)[-c^{2}_{3}(1+\frac{1}{\gamma^{2}(\left|\bar{u}_{n_{j}}\right|)})^{2}-c_{5}+\frac{k}{\gamma^{2}(\left|\bar{u}_{n_{j}}\right|)}\int^{T}_{0}H(kt,\bar{u}_{n_{j}})dt]$$
which by assumption $(H_{2})(i)$ implies that $\Phi_{k}(u_{n_{j}})\longrightarrow\infty$ as $j\longrightarrow\infty$. This contradicts the boundedness of $(\Phi_{k}(u_{n_{j}}))$. So $(\tilde{u}_{n_{j}})$ is bounded.\\
Assume that $(\bar{u}_{n_{j}})$ is unbounded, then up to a subsequence, if necessary, we can assume that $\left|\bar{u}_{n_{j}}\right|\longrightarrow\infty$ as $j\longrightarrow\infty$. As in (3.8), there exists a positive constant $c_{7}$ such
that for $j$ large enough
$$k\left|\int^{T}_{0}[H(kt,u_{n_{j}})-H(t,\bar{u}_{n_{j}})]dt\right|\leq c_{7}\gamma(\left|\bar{u}_{n_{j}}\right|).\leqno(3.9)$$
So by (3.9), we get for a positive constant $c_{8}$
$$\Phi_{k}(u_{n_{j}})\geq - c_{8}\gamma(\left|\bar{u}_{n_{j}}\right|)+k\int^{T}_{0}H(kt,\bar{u}_{n_{j}})dt$$
$$=\gamma^{2}(\left|\bar{u}_{n_{j}}\right|)[-\frac{c_{8}}{\gamma(\left|\bar{u}_{n_{j}}\right|)}+\frac{k}{\gamma^{2}(\left|\bar{u}_{n_{j}}\right|)}\int^{T}_{0}H(kt,\bar{u}_{n_{j}})dt]$$
which by assumption $(H_{2})(i)$ implies that $\Phi_{k}(u_{n_{j}})\longrightarrow\infty$ as $j\longrightarrow\infty$. This contradicts the boundedness of $(\Phi_{k}(u_{n_{j}}))$. Then $(\bar{u}_{n_{j}})$ is also bounded and therefore $(u_{n_{j}})$ is bounded. Going if necessary to a subsequence, we can assume that $\tilde{u}_{n_{j}}\rightharpoonup\tilde{u}$, $\bar{u}_{n_{j}}\longrightarrow\bar{u}$. Notice that
$$\left\|u^{+}_{n_{j}}-u^{+}\right\|=\Big(\Phi'_{k,n_{j}}(u_{n_{j}})-\Phi_{k}'(u)\Big)(u^{+}_{n_{j}}-u^{+})$$
$$-k\int^{T}_{0}(H'(kt,u_{n_{j}})-H'(kt,u)).(u^{+}_{n_{j}}-u^{+})dt$$
which implies that $u^{+}_{n_{j}}\longrightarrow u^{+}$ in $E$. Similarly, $u^{-}_{n_{j}}\longrightarrow u^{-}$ in $E$. It follows that $u_{n{j}}\longrightarrow u$ in $E$ as $j\longrightarrow\infty$ and $\Phi_{k}'(u)=0$. So $\Phi_{k}$ satisfies the $(PS)^{*}_{c}$ condition for all level $c\in\rb$. The proof of Lemma 3.1 is complete.\\
\end{Prf}
Now, let $u=u^{+}+\bar{u}\in Z$, then as in (3.8), we have for a positive constant $c_{9}$
$$k\left|\int^{T}_{0}[H(kt,u)-H(kt,\bar{u})]dt\right|\leq c_{9}\left\|u^{+}\right\|\big[ \left\|u^{+}\right\|^{\alpha}+ \gamma(\left|\bar{u}\right|)+1\big].$$
So we have
$$\Phi_{k}(u)\geq \left\|u^{+}\right\|^{2}-c_{9}\left\|u^{+}\right\|\big[ \left\|u^{+}\right\|^{\alpha}+ \gamma(\left|\bar{u}\right|)+1\big]+k\int^{T}_{0}H(kt,\bar{u})dt. \leqno(3.10)$$
Let $0<\epsilon<1$, we have
$$c_{9}\left\|u^{+}\right\|\gamma(\left|\bar{u}\right|)\leq \frac{c^{2}_{9}\gamma^{2}(\left|\bar{u}\right|)}{\epsilon^{2}}+\epsilon^{2}\left\|u^{+}\right\|^{2}. \leqno(3.11)$$
By combining (3.10) and (3.11) we get
$$\Phi_{k}(u)\geq (1-\epsilon^{2})\left\|u^{+}\right\|^{2}-c_{9}\left\|u^{+}\right\|^{1+\alpha}-c_{9}\left\|u^{+}\right\|$$
$$+\gamma^{2}(\bar{u})\big[-\frac{c^{2}_{9}}{\epsilon^{2}}+\frac{k}{\gamma^{2}(\bar{u})}\int^{T}_{0}H(kt,\bar{u})dt\big].$$
Since $0\leq\alpha<1$, we deduce by $(H_{2})(i)$ that
$$\Phi_{k}(u)\longrightarrow+\infty\ as\ \left\|u\right\|\longrightarrow\infty,\ u\in Z. \leqno(3.12)$$
Let $u\in W$ and $\xi\in \rb^{2N}$ be such that $\left|\xi\right|>0$, we have by the Mean Value Theorem, H$\ddot{o}lder's$ inequality, assumption $(H_{1})$ and the nondecreasing and properties (i), (ii) of $\gamma$
$$\left|\int^{T}_{0}[H(kt,u)-H(kt,\xi)]dt\right|=\left|\int^{T}_{0}\int^{1}_{0}H'(kt,\xi+s(u-\xi)).(u-\xi)dsdt\right|$$
$$\leq\left\|u-\xi\right\|_{L^{2}}\int^{1}_{0}\big(\int^{T}_{0}\left| H'(kt,(1-s)\xi+su)\right|^{2}dt\big)^{\frac{1}{2}}ds$$
$$\leq\left\|u-\xi\right\|_{L^{2}}\int^{1}_{0}\big(\int^{T}_{0}[p(kt)\gamma(\left|(1-s)\xi+su\right|)+q(kt)]^{2}dt\big)^{\frac{1}{2}}ds$$
$$\leq\left\|u-\xi\right\|_{L^{2}}\int^{1}_{0}\big(\int^{T}_{0}p^{2}(kt)\gamma^{2}(\left|(1-s)\xi+su\right|)dt)^{\frac{1}{2}}+\left\|q\right\|_{L^{2}}\big)ds$$
$$\leq\left\|u-\xi\right\|_{L^{2}}\Big[c(\int^{T}_{0}p^{2}(kt)[\gamma(\left|u\right|)+\gamma(\left|\xi\right|)]^{2}dt)^{\frac{1}{2}}+\left\|q\right\|_{L^{2}}\Big]$$
$$\leq\left\|u-\xi\right\|_{L^{2}}\Big[c(\int^{T}_{0}p^{2}(kt)\gamma^{2}(\left|u\right|)dt)^{\frac{1}{2}}+c\left\|p\right\|_{L^{2}}\gamma(\left|\xi\right|)+\left\|q\right\|_{L^{2}}\Big]$$
$$\leq\left\|u-\xi\right\|_{L^{2}}\Big[c\int^{T}_{0}p^{2}(kt)[a\left|u\right|^{\alpha}+b]^{2}dt)^{\frac{1}{2}}+c\left\|p\right\|_{L^{2}}\gamma(\left|\xi\right|)+\left\|q\right\|_{L^{2}}\Big]$$
$$\leq\left\|u-\xi\right\|_{L^{2}}\Big[ca(\int^{T}_{0}p^{2}(t)\left|u\right|^{2\alpha}dt)^{\frac{1}{2}}+cb\left\|p\right\|_{L^{2}}+c\left\|p\right\|_{L^{2}}\gamma(\left|\xi\right|)+\left\|q\right\|_{L^{2}}\Big]$$
$$\leq\left\|u-\xi\right\|_{L^{2}}\Big[ca\left\|p\right\|_{L^{\frac{2}{1-\alpha}}}          \left\|u\right\|^{\alpha}_{L^{2}}+cb\left\|p\right\|_{L^{2}}+c\left\|p\right\|_{L^{2}}\gamma(\left|\xi\right|)+\left\|q\right\|_{L^{2}}\Big].$$
So, by (2.1), for $\xi$ fixed there exists a positive constant $c_{10}$ such that
$$k\left|\int^{T}_{0}[H(kt,u)-H(kt,\xi)]dt\right|\leq c_{10}(\left\|u\right\|+1)(\left\|u\right\|^{\alpha}+1).$$
Therefore
$$\Phi_{k}(u)=-\left\|u\right\|^{2}+k\int^{T}_{0}[H(kt,u)-H(kt,\xi)]dt+k\int^{T}_{0}H(kt,\xi)dt$$
$$\leq -\left\|u\right\|^{2}+c_{10}(\left\|u\right\|+1)(\left\|u\right\|^{\alpha}+1)+k\int^{T}_{0}H(kt,\xi)dt.$$
Since  $0\leq\alpha<1$, then
$$\Phi_{k}(u)\longrightarrow-\infty\ as\ \left\|u\right\|\longrightarrow\infty,\ u\in W. \leqno(3.13)$$
Combining Lemma 3.1 and properties (3.12), (3.13) we deduce that the functional $\Phi_{k}$ satisfies all the assumptions of Theorem 2.1. Hence the Hamiltonian system $(\mathcal{H}_{k})$ possesses at least one $T-$ periodic solution $u_{k}$ which is a critical point of $\Phi_{k}$ and by remark 2.1, it satisfies
$$\Phi_{k}(u_{k})=C_{k}\geq\inf_{u\in Z}\Phi_{k}(\sqrt{k}e+u) \leqno(3.14)$$
where $e(t)=\frac{1}{\sqrt{\pi}}exp(\frac{2\pi}{T}tJ)e_{1}\in W$, $e_{1}$ denotes the first element of the standard basis of $\rb^{2N}$, with $x_{k}(t)=u_{k}(\frac{t}{k})$ is a $kT-$ periodic solution of $(\mathcal H)$. We will prove that the sequence $(u_{k})_{k\geq 1}$ has the following property:
$$\lim_{k\longrightarrow\infty}\frac{1}{k}\Phi_{k}(u_{k})=+\infty.\leqno(3.15)$$
This will be done by the use of some estimates on the levels $C_{k}$ of $\Phi_{k}$. For this aim the following two lemmas will be needed.\\
{\bf Lemma 3.2.}[13] Let $F:\rb\times\rb^{2N}\longrightarrow\rb$ be a continuous function $T-$ periodic in $t$ and let $C$ be a subset of $[0,T]$ with $meas(C)>0$. Assume that there exists a $T-$ periodic function $f\in L^{1}(0,T;\rb)$ such that 
$$F(t,x)\longrightarrow+\infty\ as\ \left|x\right|\longrightarrow\infty,\ a.e.\ t\in C,$$
and
$$F(t,x)\geq f(t),\ \forall x\in\rb^{2N},\ a.e.\ t\in[0,T].$$
Then for every $\delta > 0$, there exists a measurable subset $C_{\delta}$ of $C$ with $meas(C-C_{\delta})< \delta$ such that
$$F(t,x)\longrightarrow+\infty\ as\ \left|x\right|\longrightarrow\infty,\ uniformly\ in\ t\in C_{\delta}. \leqno(3.16)$$
{\bf Lemma 3.3.} Assume that $H$ satisfies $(H_{3})(i)$, then
$$\lim_{k\longrightarrow\infty}\inf_{u\in Z}\frac{\Phi_{k}(\sqrt{k}e+u)}{k}=+\infty. \leqno(3.17)$$
\begin{Prf}
Arguing by contradiction and assume that there exist sequences $k_{j}\longrightarrow\infty$, $(u_{j})\subset Z$ and a constant $c_{11}$ such that
$$\Phi_{k_{j}}(\sqrt{k_{j}}e+u_{j})\leq k_{j}c_{11},\ \forall j\in\nb. \leqno(3.18)$$
Taking $u_{j}=\sqrt{k_{j}}(u^{+}_{j}+\bar{u}_{j})$, with $u^{+}_{j}\in E^{+}$, $\bar{u}_{j}\in \rb^{2N}$, we obtain by an easy calculation
$$\Phi_{k_{j}}(\sqrt{k_{j}}e+u_{j})=k_{j}\Big[\left\|u^{+}_{j}\right\|^{2}-1+\int^{T}_{0}H(k_{j}t,\sqrt{k_{j}}(e+u^{+}_{j}+\bar{u}_{j}))dt\Big]. \leqno(3.19)$$
On the other hand, by $(H_{3})(i)$ we have
$$\int^{T}_{0}H(k_{j}t,\sqrt{k_{j}}(e+u^{+}_{j}+\bar{u}_{j}))dt\geq\int^{T}_{0}f(k_{j}t)dt=\int^{T}_{0}f(t)dt \leqno(3.20)$$
so there exists a positive constant $c_{12}$ such that
$$\Phi_{k_{j}}(\sqrt{k_{j}}e+u_{j})\geq k_{j}(\left\|u^{+}_{j}\right\|^{2}-c_{12}). \leqno(3.21)$$
The inequalities (3.18) and (3.21) imply that $(u^{+}_{j})$ is a bounded sequence in $E$. Up to a subsequence, if necessary, we can find $u^{+}\in E^{+}$ such that
$$u^{+}_{j}(t)\longrightarrow u^{+}(t)\ as\ j\longrightarrow\infty,\ a.e.\ t\in[0,T]. \leqno(3.22)$$
We claim that $(\bar{u}_{j})$ is also bounded in $E$. Indeed, if we assume otherwise, then by taking a subsequence if necessary, (3.22) implies that
$$\sqrt{k_{j}}\left|e(t)+u^{+}_{j}(t)+\bar{u}_{j}\right|\longrightarrow\infty\ as\ j\longrightarrow\infty,\ a.e.\ t\in [0,T]. \leqno(3.23)$$
Let $\delta=\frac{1}{2}meas(C)$ and $C_{\delta}$ be as defined in Lemma 3.2. For all positive integer $j$, let us define the subset $C^{j}_{\delta}$ of $[0,T]$ by
$$C^{j}_{\delta}=\frac{1}{k_{j}}\cup^{k_{j}-1}_{r=0}(C_{\delta}+rT).$$
It is easy to verify that $meas(C^{j}_{\delta})=meas(C_{\delta})$ and
$$H(k_{j}t,x)\longrightarrow+\infty\ as\ \left|x\right|\longrightarrow\infty,\ uniformly\ in\ t\in C^{j}_{\delta}. \leqno(3.24)$$
By $(H_{3})(i)$, we have
$$\int^{T}_{0}H(k_{j}t,\sqrt{k_{j}}(e+u^{+}_{j}+\bar{u}_{j}))dt$$
$$\geq\int_{C^{j}_{\delta}}H(k_{j}t,\sqrt{k_{j}}(e+u^{+}_{j}+\bar{u}_{j}))dt+\int_{[0,T]-C^{j}_{\delta}}f(k_{j}t)dt$$
$$\geq \int^{T}_{0}\chi_{C^{j}_{\delta}}H(k_{j}t, \sqrt{k_{j}}(e+u^{+}_{j}+\bar{u}_{j}))dt-\int^{T}_{0}\left|f(t)\right|dt. \leqno(3.25)$$
On the other hand, by (3.24) and Fatou's lemma, we get
$$\int_{C^{j}_{\delta}}H(k_{j}t,\sqrt{k_{j}}(e+u^{+}_{j}+\bar{u}_{j}))dt\longrightarrow +\infty\ as\ j \longrightarrow\infty\leqno(3.26)$$
so we deduce from (3.19), (3.25) and (3.26) that
$$\frac{\Phi_{k_{j}}(\sqrt{k_{j}}e+u_{j})}{k_{j}}\longrightarrow +\infty\ as\ j\longrightarrow\infty\leqno(3.27)$$
which contradicts (3.18) and proves our claim. Hence, by taking a subsequence, if necessary, we can assume that there exists $\bar{u}\in E^{0}$ such that
$$e(t)+u^{+}_{j}(t)+\bar{u}_{j}\longrightarrow u(t)=e(t)+u^{+}(t)+\bar{u}\ as\ j\longrightarrow\infty,\ a.e.\ t\in [0,T].$$
By Fourier analysis, we have $e(t)+u^{+}(t)+\bar{u}\neq 0$ for almost every $t\in [0,T]$. Therefore
$$\sqrt{k_{j}}\left|e(t)+u^{+}_{j}(t)+\bar{u}_{j}\right|\longrightarrow\infty\ as\ j\longrightarrow\infty,\ a.e.\ t\in [0,T] \leqno(3.28)$$
and by using (3.24) and Fatou's lemma, we obtain (3.27) as above, which contradicts (3.18). This concludes the proof of Lemma 3.3. \\
\end{Prf}
Now, by (3.14) and Lemma 3.3, we have
$$\frac{C_{k}}{k}\longrightarrow\infty\ as\ k\longrightarrow\infty. \leqno(3.29)$$
We claim that $\left\|u_{k}\right\|_{\infty}=\left\|x_{k}\right\|_{\infty}\longrightarrow \infty$ as $k\longrightarrow\infty$. Indeed, if we suppose otherwise, $(u_{k})$ possesses a bounded subsequence $(u_{k_{j}})$. Since
$$\frac{\Phi_{k_{j}}(u_{k_{j}})}{k_{j}}=-\frac{1}{2}\int^{T}_{0}H'(k_{j}t,u_{k_{j}}).u_{k_{j}}dt+\int^{T}_{0}H(k_{j}t,u_{k_{j}})dt$$
the sequence $(\frac{C_{k_{j}}}{k_{j}})$ is bounded, which contradicts (3.29). Consequently, we have 
$\left\|u_{k}\right\|_{\infty}\longrightarrow\infty\ as\ k\longrightarrow\infty$, which completes the proof of Theorem 3.1.\\
{\bf Proof of Corollary 3.1.} Using Lemma 3.2, it is easy to see that assumption $(H_{4})$ implies assumptions $(H
_{2})$ and $(H_{3})$. Then Corollary 3.1 is a particular case of Theorem 3.1.
\section{Minimal periods} In this section, we consider a continuous increasing function $\gamma:\rb^{+}\longrightarrow\rb^{+}$ satisfying the properties (i)-(iii) in section 3 and the following property:\\
(iv) There exists a positive constant $c_{0}$ such that for all constant $c$
$$\lim_{r\longrightarrow\infty}\frac{1}{\gamma^{2}(r)}[\int^{r}_{1}\frac{\gamma^{2}(u)}{u}du-c\ lnr]\geq c_{0}.$$
Consider the assumptions:\\
$(H'_{1})$ There exist two $T-$ periodic functions $p\in L^{\infty}(0,T;\rb^{+})$ and $q\in L^{2}(0,T;\rb^{+})$ such that
$$\left|H'(t,x)\right|\leq p(t)\gamma(\left|x\right|)+q(t),\ \forall x\in\rb^{2N},\ a.e.\ t\in [0,T].$$
$(H_{5})$ There exists a $T-$ periodic function $f\in L^{1}(0,T;\rb)$ such that either
$$\frac{H'(t,x).x}{\gamma^{2}(\left|x\right|)}\longrightarrow+\infty\ as\ \left|x\right|\longrightarrow +\infty,\ uniformly\ in\ t\in [0,T], \leqno(i)$$
or
$$\frac{H'(t,x).x}{\gamma^{2}(\left|x\right|)}\longrightarrow-\infty\ as\ \left|x\right|\longrightarrow +\infty,\ uniformly\ in\ t\in [0,T]. \leqno(ii)$$
Our main result in this section is:\\
{\bf Theorem 4.1.} Suppose $(H'_{1})$ and $(H_{5})$ hold. Then, for all positive integer $k$, the Hamiltonian system $(\mathcal H)$ possesses at least one $kT-$ periodic solution $x_{k}$ such that
$$\lim_{k\longrightarrow\infty}\left\|x_{k}\right\|_{\infty}=+\infty.$$
If moreover $H$ satisfies the following assumption
$$If\ u(t)\ is\ a\ periodic\ function\ with\ minimal\ period\ rT,$$
$$r\ rational,\ and\ H'(t,u(t))\ is\ a\ periodic\ function\ with \leqno(H)$$
$$minimal\ period\ rT,\ then\ r\ is\ necessary\ an\ integer,$$
then, for any sufficiently large prime number $k$, $kT$ is the minimal period of $x_{k}$.\\
{\bf Example 4.1.} The function $\gamma:\rb^{+}\longrightarrow\rb^{+}, t\longmapsto ln^{\frac{1}{2}}(1+t^{2})$ is a continuous nondecreasing function satisfying conditions $(i)-(iv)$. Take
$$H(t,x)=(\frac{3}{2}+sin(\frac{2\pi}{T}t))ln^{\frac{5}{2}}(1+\left|x\right|^{2}).$$
It is easy to verify that $H$ satisfies $(H'_{1})$, $(H_{5})$ and does not satisfy the assumptions of Theorem 1.2 in $[1]$. Theorem 1.2 in [1] is then a particular case of Theorem 4.1 with control function $\gamma(t)=t^{\alpha}$, $\alpha\in [0,1[$,  $t\in\rb^{+}$.\\
{\bf Proof of Theorem 4.1.} By Remark 3.1, it suffices to prove the case when $H$ satisfies $(H_{5})(i)$.\\
Firstly, let us prove that assumptions $(H'_{1})$, $(H_{5})(i)$ imply $(H_{2})(i)$. Indeed, by $(H_{5})(i)$, for all $\rho>0$, there exists a constant $c_{\rho}$ such that
$$H'(t,x).x\geq \rho\gamma^{2}(\left|x\right|)-c_{\rho},\ \forall x\in\rb^{2N},\ \forall t\in [0,T]. \leqno(4.1)$$
Let $x\in \rb^{2N}$ be such that $\left|x\right|\geq 1$, we have by the Mean Value Theorem
$$H(t,x)=H(t,0)+\int^{\frac{1}{\left|x\right|}}_{0}H'(t,sx).xds+\int^{1}_{\frac{1}{\left|x\right|}}H'(t,sx).xds.$$
By $(H'_{1})$, we have for $a.e. t\in [0,T]$
$$\left|\int^{\frac{1}{\left|x\right|}}_{0}H'(t,sx).xds\right|\leq \left|x\right|\int^{\frac{1}{\left|x\right|}}_{0}[p(t)\gamma(s\left|x\right|)+q(t)]ds\leq p(t)\gamma(1)+q(t).\leqno(4.2)$$
By (4.1), we get
$$\int^{1}_{\frac{1}{\left|x\right|}}H'(t,sx).xds=\int^{1}_{\frac{1}{\left|x\right|}}H'(t,sx).sx\frac{ds}{s} \geq\int^{1}_{\frac{1}{\left|x\right|}}[\rho\gamma^{2}(s\left|x\right|)-c_{\rho}]\frac{ds}{s}$$
$$=\rho\int^{1}_{\frac{1}{\left|x\right|}}\gamma^{2}(s\left|x\right|)\frac{ds}{s}-c_{\rho}ln(\left|x\right|)=\rho\int^{\left|x\right|}_{1}\gamma^{2}(u)\frac{du}{u}-c_{\rho}ln(\left|x\right|).\leqno(4.3)$$
Combining (4.2), (4.3), yields
$$H(t,x)\geq -p(t)\gamma(1)-q(t)+\rho\int^{\left|x\right|}_{1}\gamma^{2}(u)\frac{du}{u}-c_{\rho}ln(\left|x\right|)+H(t,0),\leqno(4.4)$$
which by property $(iv)$ of $\gamma$, imply that $H(t,x)\longrightarrow +\infty\ as\ \left|x\right|\longrightarrow\infty,\ a.e.\ t\in [0,T]$ and assumption $(H_{3})(i)$ is satisfied.
By integrating $(4.4)$, we obtain
$$\int^{T}_{0}H(t,x)dt\geq -\gamma(1)\int^{T}_{0}p(t)dt-\int^{T}_{0}q(t)dt$$
$$+T\rho \int^{\left|x\right|}_{1}\gamma^{2}(u)\frac{du}{u}-T c_{\rho}ln(\left|x\right|)+\int^{T}_{0}H(t,0)dt$$
and by property $(iv)$ of $\gamma$, we get
$$\lim_{\left|x\right|\longrightarrow\infty}\frac{1}{\gamma^{2}(\left|x\right|)}\int^{T}_{0}H(t,x)dt\geq c_{0}\rho T. \leqno(4.5)$$
Since $\rho$ is arbitrary chosen, then $H$ satisfies $(H_{2})(i)$.\\
We deduce from Theorem 3.1, that for all positive integer $k$, the system $(\mathcal H)$ possesses at least one $kT-$ periodic solution $x_{k}$ satisfying $\lim_{k\longrightarrow\infty}\left\|x_{k}\right\|_{\infty}=\infty$. It remains to study the minimality of periods of $x_{k}$, $k\geq 1$. Consider the family of functionals
$$\psi_{k}(x)=\frac{1}{2}\int^{kT}_{0}J\dot{x}.xdt+\int^{kT}_{0}H(t,x)dt$$
defined respectively on the spaces $E_{k}=H^{\frac{1}{2}}(S^{1}_{k},\rb^{2N})$ with $S^{1}_{k}=\rb/(kT\zb)$. It is easy to see that for all $k\geq 1$, $x_{k}$ is a critical point of $\psi_{k}$ and by (3.29), we have
$$\lim_{k\longrightarrow\infty}\frac{\psi_{k}(x_{k})}{k}=+\infty,. \leqno(4.6)$$
Now, let us denote by $S_{T}$ the set of $T-$periodic solutions of $(\mathcal H)$. We claim that $S_{T}$ is bounded in $H^{\frac{1}{2}}_{T}$. Indeed, assume by contradiction that there exists a sequence $(x_{n})\subset S_{T}$ such that $\left\|x_{n}\right\|\longrightarrow\infty$ as $n\longrightarrow\infty$. Let $x_{n}=x^{+}_{n}+x^{-}_{n}+x^{0}_{n}$. Multiplying both sides of the identity
$$J\dot{x}_{n}+H'(t,x_{n}(t))=0 \leqno(4.7)$$
by $x^{+}_{n}$ and integrating, we obtain
$$2\left\|x^{+}_{n}\right\|^{2}+\int^{T}_{0}H'(t,x_{n}(t))\cdot x^{+}_{n}dt=0.\leqno(4.8)$$
Using H$\ddot{o}$lder's inequality, assumption $(H'_{1})$, property $(ii)$ of $\gamma$ and inequality (2.1), we can find as above a positive constant $c_{13}$ such that
$$\left\|x^{+}_{n}\right\|\leq c_{13}\Big(\left\|x_{n}\right\|^{\alpha}+1\Big).$$
Since $0\leq\alpha <1$, this yields
$$\frac{\left\|x^{+}_{n}\right\|}{\left\|x_{n}\right\|}\longrightarrow 0\ as\ n\longrightarrow\infty. \leqno(4.9)$$
Similarly, we have
$$\frac{\left\|x^{-}_{n}\right\|}{\left\|x_{n}\right\|}\longrightarrow 0\ as\ n\longrightarrow\infty. \leqno(4.10)$$
Taking $y_{n}=\frac{x_{n}}{\left\|x_{n}\right\|}$ and using (4.9) and (4.10), we may assume without loss of generality that $y_{n}\longrightarrow y_{0}\in E^{0}$, with $\left|y_{0}\right|=1$. Since the embedding $E\longrightarrow L^{2}$, $u\longmapsto u$ is compact, we can assume, by taking a subsequence if necessary that
$$y_{n}(t)\longrightarrow y_{0}\ as\ n\longrightarrow\infty,\ a.e.\ t\in [0,T], \leqno(4.11)$$
and consequently
$$\left|x_{n}(t)\right|\longrightarrow +\infty\ as\ n\longrightarrow\infty,\ a.e.\ t\in [0,T]. \leqno(4.12)$$
So by Fatou's lemma and property (iii) of $\gamma$, we obtain
$$\int^{T}_{0}\gamma^{2}(\left|x_{n}\right|)dt\longrightarrow\infty\ as\ n\longrightarrow\infty.\leqno(4.13)$$
On the other hand, by (4.1), we have
$$\rho\int^{T}_{0}\gamma^{2}(\left|x_{n}(t)\right|)dt-c_{\rho}\leq \int^{T}_{0}H'(t,x_{n}(t)).x_{n}(t)dt. \leqno(4.14)$$
Furthermore, by Proposition 3.2 in [8], we have
$$\int^{T}_{0}H'(t,x_{n}(t))x_{n}(t)dt\leq \frac{T}{2\pi}\int^{T}_{0}\left|H'(t,x_{n}(t))\right|^{2}dt. \leqno(4.15)$$
Combining (4.14), (4.15) yields
$$\rho\int^{T}_{0}\gamma^{2}(\left|x_{n}\right|)dt-c_{\rho}\leq\frac{T}{2\pi}
\int^{T}_{0}\left|H'(t,x_{n})\right|^{2}dt.\leqno(4.16)$$
Using $(H'_{1})$ and H$\ddot{o}$lder's inequality we obtain for a positive constant $c_{14}$
$$(\int^{T}_{0}\left|H'(t,x_{n})\right|^{2}dt)^{\frac{1}{2}}\leq (\int^{T}_{0}[p(t)\gamma(\left|x_{n}\right|)+q(t)]^{2}dt)^{\frac{1}{2}}$$
$$\leq\left\|p\right\|_{\infty}\sqrt{T}(\int^{T}_{0}\gamma^{2}(\left|x_{n}\right|)dt)^{\frac{1}{2}}+\left\|q\right\|_{L^{2}}$$
$$\leq c^{\frac{1}{2}}_{14}\frac{\sqrt{2\pi}}{\sqrt{T}}\Big(\int^{T}_{0}\gamma^{2}(\left|x_{n}\right|)dt+ 1\Big)^{\frac{1}{2}}. \leqno(4.17)$$
Combining (4.16) and (4.17) yields
$$\rho\int^{T}_{0}\gamma^{2}(\left|x_{n}\right|)dt-c_{\rho}\leq c_{14}\Big(\int^{T}_{0}\gamma^{2}(\left|x_{n}\right|)dt+ 1\Big). \leqno(4.18)$$
Since $\rho > 0$ is arbitrary chosen then $(\int^{T}_{0}\gamma^{2}(\left|x_{n}\right|)dt)$ must be bounded, which contradicts (4.13). Hence $S_{T}$ is bounded and as a consequence $\psi_{1}(S_{T})$ is bounded. On the other hand, for any $x\in S_{T}$ one has $\psi_{k}(x)=k\psi_{1}(x)$, so there exists a positive constant $M$ such that
$$\forall x\in S_{T},\ \forall k\geq 1,\ \frac{\left|\psi_{k}(x)\right|}{k}\leq M. \leqno(4.19)$$
Consequently (4.6) and (4.19) show that for for all integer $k$ sufficiently large, $x_{k}\notin S_{T}$. So , by assumption $(H)$, if $k$ is chosen to be prime number, the minimal period of $x_{k}$ has to be $kT$. The proof of Theorem 4.1 is complete.\\
{\bf References} \par\noindent
[1] A. Daouas, M. Timoumi, "Subharmonics for not uniformly coercive Hamiltonian systems", Nonlinear Analysis 66 (2007) pp 571-581. \par\noindent
[2] I. Ekeland, H. Hofer, "Subharmonics for convex nonautonomous Hamiltonian systems", Comm. Pure Appl. Math. 40 (1987) pp 1-36. \par\noindent
[3] A. Fonda, A.C. Lazer, "Subharmonic solutions of conservative systems with nonconvex potentials", Proc. Amer. Math. Soc. 115 (1992) pp 183-190. \par\noindent
[4] A. Fonda, M. Ramos, M. Willem, "Subharmonic solutions for second-order differential equations", Topol. Meth. Nonl. Anal. 1 (1993) pp 49-66. \par\noindent
[5] G. Fournier, T. Timoumi, M. Willem, "The limiting case for strongly indefinite functionals", Top. Meth. in Nonlinear Analysis 1, (1993) pp 203-209.\par\noindent
[6] C.G. Liu, "Subharmonic solutions of Hamiltonian systems", Nonlinear Analysis 42 (2000) pp 185-198. \par\noindent
[7] C. Li, Z.Q. Ou, C.L. Tang, "Periodic and subharmonic solutions for a class of non-autonomous Hamiltonian systems", Nonlinear Analysis 75 (2012) pp 2262-2272.\par\noindent
[8] J. Mawhin, M. Willem, "Critical point theory and Hamiltonian systems", Springer, (1989)\par\noindent
[9] R. Michalek, G. Tarantello, "Subharmonic solutions with prescribed minimal period for nonautonomous Hamiltonian systems", J. Diff. Eq. 72 (1988) pp 28-55. \par\noindent
[10] Z.Q. Ou, C.L. Tang, "Periodic and subharmonic solutions for a class of superquadratic Hamiltonian systems", Nonlinear Analysis 58 (2004) pp 245-258. \par\noindent
[11] P.H. Rabinowitz, "Minimax methods in critical point theory with applications to differential equations",
CBMS. Reg. Conf. Ser. Math., vol 65, Amer. Math. Soc., Providence, RI (1986). \par\noindent
[12] E.A. de B.e. Silva, "Subharmonic solutions for subquadratic Hamiltonian systems", J.diff.eq. 115 (1995) pp 120-145. \par\noindent
[13] C.L. Tang, X.P. Wu, "Periodic solutions for second order systems with not uniformly coercive potential", J. Math. Anal. Appl. 259 (2001) pp 386-397. \par\noindent
[14] C.L. Tang, X.P. Wu, "Subharmonic solutions for nonautonomous sublinear second order Hamiltonian systems", J. Math. Anal. Appl. 304 (2005) 383-393. \par\noindent
[15] M. Timoumi, "Subharmonics of a Hamiltonian systems class", Dem. Math. vol XXXVII N0 4 (2004) pp 977-990. \par\noindent
[16] M. Timoumi, "Subharmonic of nonconvex Hamiltonian systems", Arch. der Math. 73 (1999) pp 422-429. \par\noindent
[17] M. Timoumi, "Subharmonic oscillations of a class of Hamiltonian systems", Nonlinear analysis 68, (2008) pp 2697-2708. \par\noindent
[18] M. Timoumi, "Subharmonic solutions for nonautonomous second order Hamiltonian systems", Elect. J. Diff. Eq., Vol 2012, (2012).

\end{document}